# Hölder Continuity and Differentiability Almost Everywhere of $(K_1, K_2)$-Quasiregular Mappings

GAO Hongya[1]    LIU Chao[1]    LI Junwei[2,1]

1. College of Mathematics and Computer Science, Hebei University, Baoding, 071002, China
2. Information Center, Hebei Normal College for Nationalities, Chengde, 067000, China

**Abstract.** This paper deals with $(K_1, K_2)$-quasiregular mappings. It is shown, by Morrey's Lemma and isoperimetric inequality, that every $(K_1, K_2)$-quasiregular mapping satisfies a Hölder condition with exponent $\alpha$ on compact subsets of its domain, where

$$\alpha = \begin{cases} 1/K_1, & \text{for } K_1 > 1, \\ \text{any positive number less than } 1, & \text{for } K_1 = 1 \text{ and } K_2 > 0, \\ 1, & \text{for } K_1 = 1 \text{ and } K_2 = 0, \\ 1, & \text{for } K_1 < 1. \end{cases}$$

Differentiability almost everywhere of $(K_1, K_2)$-quasiregular mappings is also derived.

## §1  Introduction and Statement of Results

Let $\Omega$ be an arbitrary open set in $\mathrm{R}^n$, $n \geq 2$. For any point $x \in \Omega$ and $r > 0$, we denote by $B(x, r)$ the ball with radius $r$ centered at $x$ and $S(x, r) = \partial B(x, r)$ the sphere of $B(x, r)$. Let $|B(x, r)| = \omega_n r^n$ be the $n$-dimensional Lebesgue measure of the ball $B(x, r)$, where $\omega_n$ be the volume of the unit ball in $\mathrm{R}^n$. Denote by $\rho_\Omega(x) = \text{dist}(x, \partial\Omega)$ the distance from $x$ to $\partial\Omega$, with the subscript $\Omega$ omitted whenever no confusion can result.

For $p \geq 1$, we denote by $L^p(\Omega)$ the $L^p$ space of functions on $\Omega$, $W^{1,p}(\Omega)$ will denote the corresponding Sobolev space of functions in $L^p(\Omega)$ whose distributional first derivatives belong also to the space $L^p(\Omega)$. Similarly, $W^{1,p}(\Omega, \mathrm{R}^n)$ will be the space of functions $f = (f^1, f^2, \cdots, f^n) : \Omega \to \mathrm{R}^n$ such that $f^i \in W^{1,p}(\Omega)$ for $i = 1, 2, \cdots, n$.

For $A$ an $n \times n$ matrix, we define the norm of $A$ as $|A| = \sup_{|\xi|=1} |A\xi|$.

Research supported by National Natural Science Foundation of China (Grant No. 10971224) and Natural Science Foundation of Hebei Province (Grant No. A2011201011).

We say that a function $f : \Omega \to \mathrm{R}^n$ does not change sign in $\Omega$ if either $u(x) \geq 0$ almost everywhere in $\Omega$ or $u(x) \leq 0$ almost everywhere in $\Omega$.

A mapping $f : \Omega \to \mathrm{R}^n$ is said to be satisfy Hölder condition with exponent $\alpha$ on compact subsets of $\Omega$, where $0 < \alpha \leq 1$, if for every compact set $V \subset\subset \Omega$ there is a number $M(V)$, $0 \leq M(V) < +\infty$, such that for any $x_1, x_2 \in V$,

$$|f(x_1) - f(x_2)| \leq M(V)|x_1 - x_2|^\alpha.$$

If $f$ satisfies a Hölder condition with exponent $\alpha = 1$ on compact subsets of $\Omega$, then $f$ is said to satisfy a Lipschitz condition on compact subsets of $\Omega$.

Let $f = (f^1, f^2, \cdots, f^n) : \Omega \to \mathrm{R}^n$ be a mapping in $W^{1,n}(\Omega, \mathrm{R}^n)$. The linear mapping

$$Df(x) = \left(\frac{\partial f^i}{\partial x_j}\right)_{1 \leq i,j \leq n}$$

is defined for almost all $x \in \Omega$. Its determinate $\det Df(x)$ is called the Jacobian of $f$ at the point $x$, and is denoted by $\mathcal{J}(x, f)$.

In [1], Zheng and Fang gave the definition for $(K_1, K_2)$-quasiregular and quasiconformal mappings.

**Definition 1.1.** *A mapping $f = (f^1, f^2, \cdots, f^n) : \Omega \to R^n$ is called $(K_1, K_2)$-quasiregular with $0 < K_1 < +\infty$, $0 \leq K_2 < +\infty$, if it satisfies the following conditions:*
  (i)  *$f$ belongs to the class $W^{1,n}(\Omega, R^n)$,*
  (ii) *the Jacobian $\mathcal{J}(x, f)$ does not change sign in $\Omega$, and*

$$|Df(x)|^n \leq K_1|\mathcal{J}(x,f)| + K_2 \tag{1.1}$$

*for almost all $x \in \Omega$.*

*A mapping $f = (f^1, f^2, \cdots, f^n) : \Omega \to R^n$ is said to be $(K_1, K_2)$-quasiconformal if it satisfies (i), (ii), and*
  (iii) *$f$ is a homeomorphism.*

The estimate of the modulus of continuity of $(K_1, 0)$-quasiconformal mappings was first established by Kreines [2]. The Hölder property was first proved for a $(K_1, 0)$-quasiregular mapping by Reshetnyak [3,4], and simultaneously by Callender [5]. Simon [6] established an estimate of Hölder continuity when he studied $(K_1, K_2)$-quasiconformal mappings between two surfaces of the Euclidean space $\mathrm{R}^3$. This estimate has important applications to elliptic equations with two variables. In [7], Gilbarg and Trudinger obtained an *a priori* $C^{1,\alpha}_{loc}$ estimate for quasilinear elliptic equations with two variables by using the Hölder continuity method established in the study of plane

$(K_1, K_2)$-quasiregular mappings, and then established the existence theorem of Dirichlet problem. Many results on quasiregular mappings and their applications to nonlinear PDEs and elasticity theory have been established recently, see [8-10] and the references therein.

Because of the importance of plane $(K_1, K_2)$-quasiregular mappings to the *a priori* estimates in nonlinear PDE theory, Zheng and Fang [1] developed the theory of $(K_1, K_2)$-quasiregular mappings in 1998 by using the theory of outer differential forms and Grassman algebra, and obtained an $L^p$-integrability $(p > n)$ result for space $(K_1, K_2)$-quasiregular mappings. For some other developments on $(K_1, K_2)$-quasiregular mapping theory, see [11-15].

It is a typical situation in quasiconformal analysis that one wants to build up the Hölder continuity theory for $(K_1, K_2)$-quasiregular mappings. In this paper, we generalize the results of [1,11,13], and the following Hölder continuity result is obtained.

**Theorem 1.1.** *Let $f : \Omega \to R^n$ be a $(K_1, K_2)$-quasiregular mapping. Assume that*
$$\int_\Omega |Df(x)|^n dx = M < +\infty.$$
*Then the function $f$ satisfies a Hölder condition with exponent $\alpha$ on compact subsets of $\Omega$, where*
$$\alpha = \begin{cases} 1/K_1, & \text{for } K_1 > 1, \\ \text{any positive number less than } 1, & \text{for } K_1 = 1 \text{ and } K_2 > 0, \\ 1, & \text{for } K_1 = 1 \text{ and } K_2 = 0, \\ 1, & \text{for } K_1 < 1. \end{cases} \quad (1.2)$$
*Further, if $V$ is contained strictly inside $\Omega$, then for any $x, y \in V$*
$$|f(x) - f(y)| \le L|x - y|^\alpha,$$
*where the constant $L$ depends only on $V$, the constants $K_1$ and $K_2$, the dimension $n$, the distance from $V$ to the boundary of $\Omega$, and the constant $M$.*

**A counterexample** The mapping $f$ with $f(0) = 0$ and $f : x \mapsto x|x|^{\alpha - 1}$ for $x \ne 0$, where $\alpha = 1/K_1$, shows that the exponent $1/K_1$ in Theorem 1.1 is optimal. For this $f$ we have $|f(x) - f(0)| = |x|^\alpha$.

The following corollary is a direct consequence of Theorem 1.1.

**Corollary 1.1.** *Let $\Omega$ be an open domain in $R^n$, and $F(\Omega, K_1, K_2, M)$ the collection of all $(K_1, K_2)$-quasiregular mappings $f$ on $\Omega$ such that*
$$\int_\Omega |Df(x)|^n dx \le M.$$

*Then the set of functions f is equi-uniformly continuous on every compact subset of $\Omega$.*

**Definition 1.2.** *The mapping f is said to have property N if the image of every set $E \subset \Omega$ of measure zero is a set of measure zero.*

**Corollary 1.2.** *Let $f : \Omega \to R^n$ be a $(K_1, K_2)$-quasiregular mapping with $0 < K_1 < 1$, or $K_1 = 1$ and $K_2 = 0$, then f has property N.*

*Proof.* [16, Theorem 2.2] states that every locally Lipschitz mapping has property $N$, which together with Theorem 1.1 yields the desired result. □

**Definition 1.3.** *A mapping $f : \Omega \to R^n$ is said to be differentiable at a point $a \in \Omega$ if there exists a linear mapping $L : R^n \to R^n$ such that*
$$f(x) = f(a) + L(x - a) + \beta(x)|x - a|$$
*for all $x \in \Omega$, where $\beta(x) \to 0$ as $x \to a$. The mapping L is called the differential of f at the point a.*

The following theorem states that any $(K_1, K_2)$-quasiregular mapping $f$ is differential almost everywhere.

**Theorem 1.2.** *Let f be a $(K_1, K_2)$-quasiregular mapping. Then for almost all $x \in \Omega$ the linear mapping $Df(x)$ is the differential of f at the point x.*

*Proof.* The proof of Theorem 1.2 is almost line by line of the proof of [16, Theorem 1.2] by using Corollary 1.1. We omit the details. □

## §2  Preliminary Lemmas

The proof of Theorem 1.1 is based on two facts. The first is a lemma due to Morrey. The second is an isoperimetric inequality due to Reshetnyak.

**Lemma 2.1.** (Morrey's Lemma [17]) *Let $\Omega \subset R^n$ be an open subset, and $f : \Omega \to R^k$ a function of the class $W^{1,m}(\Omega, R^k)$, where $1 \leq m \leq n$. Assume that there exist numbers $\alpha(0 < \alpha \leq 1)$, $M < +\infty$, and $\delta > 0$ such that*
$$\int_{B(a,r)} |Df(x)|^m dx \leq M r^{n-m+m\alpha} \tag{2.1}$$
*for every ball $B(a,r) \subset \Omega$ with radius at most $\delta$. Then there exists a continuous function $\tilde{f}$ such that $f(x) = \tilde{f}(x)$ almost everywhere, and the oscillation of $\tilde{f}$ on any ball $B(x,r) \subset \Omega$ with $r \leq \delta/3$ and $r < \rho(x)/3$ does not exceed $CM^{1/m} r^\alpha$, where $C < +\infty$ is a constant.*

**Lemma 2.2.** (Isoperimetric Inequality [16]) *Suppose that $\Omega \subset R^n$ and the mapping $f : \Omega \to R^n$ is in the class $W^{1,n}(\Omega, R^n)$. Then for any $a \in U$ and almost all $t \in (0, \rho(a))$*

$$\int_{B(a,t)} \mathcal{J}(x,f)dx \leq \frac{t}{n}\int_{S(a,t)} |Df(x)|^n d\sigma(x), \qquad (2.2)$$

*where $d\sigma$ is the area element of the sphere $S(a,t)$.*

With the Morrey's Lemma and isoperimetric inequality in hands, we can now prove the following two lemmas, which will be used in the proof of Theorem 1.1.

**Lemma 2.3.** *Suppose that $\Omega \subset R^n$ is an open set and $f$ a $(K_1, K_2)$-quasiregular mapping. For $x \in \Omega$ and $r < \rho(x)$ let*

$$v(x, K_1, K_2, n, r) = \begin{cases} \frac{w(r)}{r^{n/K_1}} + \frac{K_2\omega_n}{K_1-1}r^{n-n/K_1}, & \text{for } K_1 \neq 1, \\ \frac{w(r)}{r^n} + K_2 n\omega_n \ln r, & \text{for } K_1 = 1, \end{cases} \qquad (2.3)$$

*where*

$$w(r) = \int_{B(x,r)} |Df(x)|^n dx.$$

*Then the function $r \mapsto v(x, K_1, K_2, n, r)$ is nondecreasing.*

*Proof.* For $r < \rho(x)$, (1.1) leads to

$$\begin{aligned}\int_{B(x,r)} |Df(x)|^n dx &\leq K_1 \int_{B(x,r)} |\mathcal{J}(x,f)|dx + K_2 \int_{B(x,r)} dx \\ &= K_1 \left|\int_{B(x,r)} \mathcal{J}(x,f)dx\right| + K_2|B(x,r)|\end{aligned} \qquad (2.4)$$

because $\mathcal{J}(x,f)$ does not change sign in $\Omega$. On the basis of Lemma 2.2

$$\left|\int_{B(x,r)} \mathcal{J}(x,f)\right| \leq \frac{r}{n}\int_{S(x,r)} |Df(x)|^n d\sigma(x) \qquad (2.5)$$

for almost all $r \in (0, \rho(x))$. From (2.4) and (2.5) we get

$$\int_{B(x,r)} |Df(x)|^n dx \leq \frac{K_1 r}{n}\int_{S(x,r)} |Df(x)|^n d\sigma(x) + K_2\omega_n r^n. \qquad (2.6)$$

Let

$$\int_{S(x,r)} |Df(x)|^n d\sigma(x) = s(r).$$

Applying Fubini's theorem, we get that $w(r) = \int_0^r s(t)dt$ for all $r \in (0, \rho(x))$. This leads us to conclude that the function $w$ is absolutely continuous and $w'(r) = s(r)$ for almost all $r \in (0, \rho(x))$. From (2.6) we have that

$$w(r) \leq \frac{K_1 r w'(r)}{n} + K_2\omega_n r^n$$

for almost all $r$. This is equivalent to
$$\frac{K_1 r w'(r)}{n} - w(r) + K_2 \omega_n r^n \geq 0.$$
Multiplying both sides of this inequality by $r^{-(n/K_1)-1}$ yields
$$\frac{K_1 w'(r)}{n r^{n/K_1}} - \frac{w(r)}{r^{n/K_1+1}} + K_2 \omega_n r^{n-n/K_1-1} \geq 0.$$
We get after obvious transformations that
$$\frac{\partial v(x, K_1, K_2, n, r)}{\partial r} \geq 0,$$
where $v(x, K_1, K_2, n, r)$ is defined by (2.3). Consequently, the function $r \mapsto v(x, K_1, K_2, n, r)$ is nondecreasing, as desired. $\square$

**Lemma 2.4.** *Suppose that $\Omega \subset R^n$ is an open set and $f$ a $(K_1, K_2)$-quasiregular mapping. Let $\int_\Omega |Df(x)|^n dx = M^n$. Then the vector-valued function $f$ is equivalent, in the sense of the theory of integral, to some continuous function $\tilde{f}$. Further, for every set $V$ lying strictly inside $\Omega$ the oscillation of $\tilde{f}$ on any ball $B(a, r)$ of radius $r < 2d/3$ about an $a \in V$ does not exceed $Cr^\alpha$, where $d = \text{dist}(V, \partial\Omega)$.*

*Proof.* Let $a \in V$ and
$$w(a, r) = \int_{B(a,r)} |Df(x)|^n dx \leq M^n.$$
According to Lemma 2.3, the function $r \mapsto v(a, K_1, K_2, n, r)$ is nondecreasing. We now divide the proof into four cases.

*Case 1* $K_1 > 1$. In this case,
$$\begin{aligned} v(a, K_1, K_2, n, r) &= \frac{w(a,r)}{r^{n/K_1}} + \frac{K_2 \omega_n}{K_1 - 1} r^{n-n/K_1} \\ &\leq v(a, K_1, K_2, n, 2d/3) \leq M^n (2d/3)^{-n/K_1} + \frac{K_2 \omega_n}{K_1 - 1}(2d/3)^{n-n/K_1} \\ &:= C_1(K_1, K_2, M, d, n). \end{aligned} \quad (2.7)$$
for all $r \in (0, 2d/3)$; hence
$$w(a, r) \leq C_1 r^{n/K_1} - \frac{K_2 \omega_n}{K_1 - 1} r^n \leq C_1 r^{n/K_1}.$$

*Case 2* $K_1 = 1$ and $K_2 > 0$. It is no loss of generality to assume that $d < 3/2$,
$$\begin{aligned} v(a, K_1, K_2, n, r) &= \frac{w(a,r)}{r^n} + K_2 n \omega_n \ln r \\ &\leq v(a, K_1, K_2, n, 2d/3) \leq M^n (2d/3)^{-n} + K_2 n \omega_n \ln(2d/3). \end{aligned} \quad (2.8)$$

This implies, for any $0 < \alpha < 1$,
$$\begin{aligned} w(a,r) &\leq [M^n(2d/3)^{-n} + K_2 n\omega_n \ln(2d/3)]r^n - K_2 n\omega_n r^n \ln r \\ &= M^n(2d/3)^{-n}r^n + K_2 n\omega_n r^n \ln(2d/3r) \\ &= \left[M^n(2d/3)^{-n}r^{n(1-\alpha)} + K_2 n\omega_n r^{n(1-\alpha)} \ln(2d/3r)\right] r^{n\alpha}. \end{aligned} \quad (2.9)$$
since $\lim_{r\to 0^+} r^{n(1-\alpha)} \ln(2d/3r) = 0$ for any $0 < \alpha < 1$, then we take $\delta$ such that $\delta^{n(1-\alpha)} \ln(2\delta/3d) \leq 1$. When $0 < r \leq \delta$, we have from (2.9) that
$$w(a,r) \leq \left[\left(\frac{3M}{2}\right)^n d^{-n\alpha} + 1\right] r^{n\alpha} := C_2 r^{n\alpha}.$$

*Case 3* $K_1 = 1$ and $K_2 = 0$. (2.8) implies
$$w(a,r) \leq M^n(2d/3)^{-n} r^n := C_3 r^n.$$

*Case 4* $K_1 < 1$. In this case, (2.7) also holds for all $r \in (0, 2d/3)$, thus
$$\begin{aligned} w(a,r) &\leq C_1 r^{n/K_1} + \frac{K_2 \omega_n}{1-K_1} r^n = \left[C_1 r^{n/K_1 - n} + \frac{K_2 \omega_n}{1-K_1}\right] r^n \\ &\leq \left[C_1 (2d/3)^{n(1-K_1)/K_1} + \frac{K_2 \omega_n}{1-K_1}\right] r^n := C_4 r^n \end{aligned}$$

In all the cases we have derived that for $0 < r \leq \delta$,
$$\int_{B(a,r)} |Df(x)|^n dx \leq Cr^\alpha,$$
where $\alpha$ is defined as (1.2) and $C$ depends only on $K_1, K_2, M, d, n$. The required result follows directly from Lemma 2.1. □

## §3  Proof of Theorem 1.1

*Proof.* Let $V$ be a compact subset of $\Omega$, and let $\gamma$ be the smaller of the numbers $\delta/3$ and $d/3$, where $\delta$ is the constant in Lemma 2.1 and $d = \text{dist}(V, \partial\Omega)$. We consider the function $h$ defined as follows on the product $V \times V$: $h(x,y) = |f(x) - f(y)|/|x-y|^\alpha$ for $x \neq y$, and $h(x,x) = 0$. Let $H$ be the set of pairs $(x,y) \in V \times V$ such that $|x-y| \geq \gamma$, and let $G = (V \times V) \setminus H$. The set $H$ is compact, and thus $h$ is bounded on $H$ by continuity. The conclusion of Lemma 2.4 enables us to deduce that $h$ is bounded also on $G$. Consequently, $h$ is bounded on $V \times V$, and thus $|f(x) - f(y)| \leq L|x-y|^\alpha$ for any $x, y \in V$. □